\documentclass[11pt,eqno]{article}
\usepackage{mathrsfs}
\usepackage{graphicx}
\usepackage{amsmath}
\usepackage{amsthm}
\usepackage{amssymb}
\usepackage{amsfonts}
\usepackage{epsfig}
\def\Im {\mathop{\rm Im}\nolimits}
\def\arg {\mathop{\rm arg}\nolimits}
\def\Re {\mathop{\rm Re}\nolimits}

\newtheorem{pro}{Proposition}
\newtheorem{cor}{Corollary}

\newtheorem{thm}{Theorem}
\numberwithin{equation}{section}

\title{{Uniform Treatment of Darboux's Method \\ and the  Heisenberg Polynomials}}

\author{{{Sai-Yu Liu}}$^{\natural}$ , {{R. Wong}}$^{\sharp}$ and {{Yu-Qiu
Zhao}}$^{\natural}$\footnote{Corresponding author.
 {\it{E-mail
address:}} {stszyq@mail.sysu.edu.cn}} \footnote{The research of Y.-Q. Zhao was
supported in part by the National Natural Science Foundation of
China under grant numbers 10471154 and 10871212.}}

\date{{\it{$^{\natural}$Department of
Mathematics, ZhongShan University, GuangZhou 510275,  China\\
$^{\sharp}$Department of Mathematics, City University of Hong Kong,
Tat Chee Avenue, Kowloon, Hong Kong}}}

\begin{document}
\maketitle

\noindent \hrule width 6.27in\vskip .3cm

\noindent {\bf{Abstract }} We show that the set of  Heisenberg
polynomials furnishes  a simple  non-trivial example in the uniform
treatment of Darboux's method. \vskip .2cm

\vskip .2cm

  \vskip .3cm

\noindent \hrule width 6.27in\vskip 1.3cm

\section*{1. The Heisenberg Polynomials} \indent\setcounter{section} {1}
The  Heisenberg polynomials are polynomials in $z$ and $\bar z$,
defined   explicitly as
\begin{equation}\label{(1.1)}
 C_n^{(\alpha, \beta)}(z) = \sum^n_{j=0}\frac{(\alpha)_j (\beta)_{n-j}}{j!(n-j)!} \bar z^j
 z^{n-j},~~~
 n=0,1,\cdots ,\end{equation}
 where $\alpha$ and $\beta$ are real constants, and $(\gamma)_k$ is  the Pochhammer
symbol  defined by  $(\gamma)_0=1$ and
$(\gamma)_k=\gamma(\gamma+1)\cdots (\gamma+k-1)$.  The
representation can be readily derived from the generating function
\begin{equation}\label{(1.2)}
    (1-w\bar{z})^{-\alpha}(1-wz)^{-\beta}=\sum_{n=0}^{\infty}C_{n}^{(\alpha,\beta)}(z)w^n,\quad\quad
    |wz|<1.
\end{equation}
The notation $C_n^{(\alpha, \beta)}(z)$ was  used by Gasper in
\cite{gasper1981}  in the sense of (\ref{(1.1)}) and (\ref{(1.2)}),
and yet the term {\it{Heisenberg polynomials}} was first used by Dunkl in
\cite{dunkl1984}.

The study of   Heisenberg polynomials may be  traced  back to the
investigation of the harmonic polynomials on Heisenberg group
$H_{N}$.  The group   $H_{N}$ consists of $(N+1)$-tuples in $\mathbb{C}^{N}\times
\mathbb{R}$, $N=1,2,3\cdots$, with the group operation
\begin{equation}\label{(1.3)}
    (\zeta,s)\cdot (z,t)=(\zeta+z,s+t+2\Im\langle\zeta,z\rangle),
\end{equation}
where $\zeta,z\in \mathbb{C}^{N}$, $s,t\in \mathbb{R}$, and
$\langle\zeta,\xi\rangle=\sum\limits_{j=1}^{N}\zeta_{j}\bar{\xi_{j}}$;
cf. Dunkl \cite{dunkl1984} and  Greiner \cite{greiner1980}.

The Heisenberg group $H_{N}$ arose  in various  branches of
mathematics, such as harmonic  analysis
\cite{howe1989,thangavelu1998} and   partial differential operators
and equations \cite{folland1973,follandstein1974}; the reason may
partially be that
the   results in
Abelian harmonic analysis could be   adapted to
Heisenberg groups.

The so-called  $L_{2\gamma}$-harmonic polynomials are the polynomial solutions
$p=p(z_1, \bar z_1,\cdots z_N, \bar z_N, t)$ of  $L_{2\gamma}p=0$ in $H_{N}$,
 where
 the differential operator  $L_{2\gamma}$ is defined by
\begin{equation}\label{(1.4)}
    L_{2\gamma}=-\frac{1}{2}\sum_{j=1}^{N}(Z_{j}\overline{Z}_{j} -\overline Z_j Z_{j})
    +2i\gamma\frac{\partial}{\partial
t},
\end{equation}
and
 $$Z_{j}=\frac{\partial}{\partial
z_{j}}+i\bar{z}_{j}\frac{\partial}{\partial t},~~~
 \overline{Z}_{j}=\frac{\partial}{\partial
\bar{z}_{j}}-i\bar{z}_{j}\frac{\partial}{\partial t},$$
  $\gamma$ being
  a constant.  In particular, $L_{0}$ in $H_1$ is a sub-elliptic Laplacian
on $\mathbb{R}^{3}$,  cf. \cite{folland1973,follandstein1974}.

Greiner \cite{greiner1980} initiated the study
 of
$L_{2\gamma}$-spherical harmonic polynomials. It is later shown by
Dunkl (cf. \cite{dunkl1984}) that every such polynomial can be
expressed in terms of  $C_{n}^{(\alpha,\beta)}$.

A few facts are known about the Heisenberg polynomials. For example,
the  representation (\ref{(1.1)}) of $C_{n}^{(\alpha,\beta)}$ can be
interpreted as a Gauss hypergeometric  function.
 The polynomials $C_{n}^{(\alpha,\beta)}(z)$ reduce to the   Legendre polynomials on the unit circle
 when   $\alpha= \beta=\frac{1}{2}$.  The readers
are referred to \cite{dunkl1984} for more properties of these
polynomials.

We mention several relevant investigations on the Heisenberg
polynomials.  In 1981, Gasper  \cite{gasper1981} indicated that
 $\{C_{n}^{(\alpha,\beta)}(e^{i\theta})\}$ furnishes  a system of orthogonal
functions on $[0,2\pi]$, with weight function
$w^{(\alpha,\beta)}(\theta)=\left (1-e^{-2i\theta}\right
)^{\alpha}\left (1-e^{2i\theta}\right )^{\beta}$, where $\alpha$,
$\beta$, $\alpha+\beta>-1$.

  In 1984,
 Dunkl  \cite{dunkl1984}    obtained a bi-orthogonality relation for these
functions in $L^{2}(T,|\sin\theta|^{\alpha+\beta}d\theta)$, where
$\alpha>0$, $\beta>0$, and $T=\{ \left. e^{i\theta}\right |
-\pi<\theta\leq \pi\}$. A 1986-paper \cite{dunkl1986} of Dunkl
addressed the density problem for
$\left\{C_{n}^{(\alpha,\beta)}(e^{i\theta})\right\}$ in $L_2$.

In 1986,  Temme   \cite{temme1986}   considered  two sets of
bi-orthogonal polynomials, $P_{n}(z;\alpha,\beta)$ and
$Q_{n}(z;\alpha,\beta)$, closely related to the Heisenberg
polynomials, where
$$P_{n}(z;\alpha,\beta)=F(-n,\alpha+\beta+1;2\alpha+1;1-z)~~\mbox{and}~~
 Q_{n}(z;\alpha,\beta)=P_{n}(z;\alpha,-\beta),$$ with $F(a,b;c;z)$
being the Gauss hypergeometric  function;  cf.
\cite[p.384]{olver2010}. The connection with the Heisenberg
polynomials is
\begin{equation}\label{(1.5)}P_{n}(e^{i\theta};\alpha,\beta)=\frac{n!}{(2\alpha+1)_{n}}e^{in\theta/2}
C_{n}^{(\alpha-\beta,\alpha+\beta+1)}(e^{i\theta/2}).\end{equation} The two sets
$\{P_{n}\}$ and $\{Q_n\}$
 are bi-orthogonal on the unit circle with respect to the weight function
$(1-e^{i\theta})^{\alpha+\beta}(1-e^{-i\theta})^{\alpha-\beta}$. Based on an integral  representation of the Gauss hypergeometric function,
uniform    asymptotic expansions for $P_n(z; \alpha, \beta)$ were   obtained in
\cite{temme1986},
  including   error bounds. These expansions are  in terms of  confluent hypergeometric functions,
 and show  that all points on the unit circle in the $z$-plane are transition points
in the sense that the polynomials have different asymptotic behavior inside and  outside   the unit disk.
To ensure the convergence of the integral representation,
  extra restrictions on the parameters, namely
$\alpha+\beta >-1$ and $\alpha-\beta \geq 0$,  are required in  \cite{temme1986}.

It is readily seen that  the Heisenberg polynomials are homogeneous
functions in  $z$ and $\bar{z}$, which implies that $C_n^{(\alpha, \beta)}(\rho e^{i\theta} )=\rho^n C_n^{(\alpha, \beta)}( e^{i\theta} )$.
 Hence, to study the behavior of the
polynomials, it is  natural to consider the polynomials on the unit
circle. The generating function (\ref{(1.2)}) now takes the form
 \begin{equation}\label{(1.6)}
   \left  (1-ze^{-i\theta}\right )^{-\alpha}\left (1-ze^{i\theta}\right )^{-\beta}
   =\sum_{n=0}^{\infty}C_{n}^{(\alpha,\beta)}(e^{i\theta})z^n
\end{equation}
for $|z|<1$ and $\theta\in [0, \pi]$.

In what follows, we will show that the set of Heisenberg polynomials
on the unit circle provides  a simple  non-trivial example in the
uniform treatment of   Darboux's method given in
\cite{wongzhao2005}.

  \noindent
 \section*{2. Uniform Treatment of Darboux's Method} \indent\setcounter{section} {2}
\setcounter{equation} {0} \label{sec:2}

Darboux was the first to consider the asymptotic behavior of the
coefficients    $a_n$ in  the Maclaurin expansion
 \begin{equation}\label{(2.1)}
   F(z)   =\sum_{n=0}^{\infty}a_n z^n,
\end{equation}where $F(w)$ has only a finite number of singularities
on its circle of convergence, all of which are algebraic in nature.
A method was introduced to obtain the asymptotic expansion for $a_n$
as $n\rightarrow \infty$. The method, named after Darboux, indicates
that the contribution to the expansion comes from the singularities
on the circle of convergence; cf., e.g., Olver \cite{olver1974} and
Wong \cite{wong1989} for the classical Darboux's method.

If the singularities are allowed  to move around on the circle of
convergence, Darboux's method will continue to work provided that  their
essential configuration remains the same (although  their relative positions
vary). However, this method breaks down, when two or more singularities coalesce
with each other. A uniform treatment was later given  by Fields
\cite{fields1968}, in which he considered   generating
functions of the type
\begin{equation}\label{(2.2)}
   F(z,\theta) =(1-z)^{-\lambda} \left [\left(e^{i\theta}-z\right )
   \left ( e^{-i\theta}-z\right )\right ]^{-\Delta} f(z, \theta)  =\sum_{n=0}^{\infty}a_n(\theta)
   z^n
\end{equation}
for $|z|<1$, where $\lambda$ and $\Delta$ are bounded quantities,
the branches are chosen such that both  $(1-z)^{-\lambda}$ and
$\left [\left(e^{i\theta}-z\right )
   \left ( e^{-i\theta}-z\right )\right ]^{-\Delta}$ reduce to $1$
   at $z=0$, and $f(z, \theta)$ is analytic in $|z|<e^{\eta}$
   with $\eta>0$   not depending  on $\theta$. The results in
   \cite{fields1968} were regarded to be too complicated to apply; see, e.g.,
  \cite[p.167]{erdelyi1970}.

The problem of uniform treatment was revisited by Wong and Zhao
\cite{wongzhao2005}, in which they firstly addressed   a special case of (\ref{(2.2)}) with $\lambda=0$, namely
\begin{equation}\label{(2.3)}
   F(z,\theta) =  \left [\left(e^{i\theta}-z\right )
   \left ( e^{-i\theta}-z\right )\right ]^{-\alpha} f(z, \theta)  =\sum_{n=0}^{\infty}a_n(\theta)
   z^n .
\end{equation}
Uniform asymptotic expansion was derived in terms of  the Bessel functions
$J_{\alpha\pm   1/ 2}(n\theta)$, valid for $\theta\in [0,
\pi-\varepsilon]$, $\varepsilon >0$. The general case
considered in \cite{wongzhao2005} deals with many coalescing
singularities
\begin{equation}\label{(2.4)}
   F(z,\theta) =  \left \{ \prod_{k=1}^q (z_k(\theta)-z)^{-\alpha_k} \right \} f(z, \theta)
   =\sum_{n=0}^{\infty}a_n(\theta)
   z^n,
\end{equation}
where $|z_k(\theta)|=1$, $q\geq 2$, $\alpha_k$ are constants,
and $f(z,\theta)$ is the same function as in (\ref{(2.2)}) and
(\ref{(2.3)}). In  (\ref{(2.4)}),    all the branch points $z_k(\theta)$ approach
$1$ as $\theta$ tends to $0$;  more precisely,
$z_k(\theta)=e^{i\theta s_k(\theta)}$ with $s_k(0)=-iz_k'(0)$ not
vanishing.  A  uniform asymptotic expansion for the general case was
presented  in \cite{wongzhao2005}, in terms of a set of special functions $T_l$,
$l=1,2,\cdots,q$.

In a follow-up paper,  Bai and Zhao \cite{baizhao2007} have extended
 the essential ideas in \cite{wongzhao2005} to   study
the classical Jacobi polynomials.

As mentioned in the last section and shown in (\ref{(1.6)}), the set
of  Heisenberg polynomials   serves as a good example of the
uniform Darboux treatment, which is simple (only two singularities) and non-trivial (the exponents $\alpha$ and
$\beta$ may be different). To
derive  the asymptotic expansion of
$C_{n}^{(\alpha,\beta)}(e^{i\theta})$, we   apply  Theorem 2 of
\cite{wongzhao2005}. To this end, we define
\begin{equation}\label{(2.5)}T_{l}(x)=\frac{1}{2\pi
i}\int_{\Gamma_{0}}s^{l-1}(s-i)^{-\beta}(s+i)^{-\alpha}e^{xs}ds,~~l=1,2, \end{equation}
where $\Gamma_{0}$ is a Hankel-type
loop which starts and ends at $-\infty$, and encircles $s=\pm i$ in
the positive sense;  cf. \cite[(4.7)]{wongzhao2005}. We further introduce an auxiliary function
\begin{equation}\label{(2.6)}h_{0}(s,\theta)=
e^{i\theta(\alpha-\beta)}\left [\frac{e^{i\theta}-e^{-\theta
s}}{(s+i)\theta}\right ]^{-\alpha}\left
[\frac{e^{-i\theta}-e^{-\theta s}}{(s-i)\theta}\right ]^{-\beta}
\end{equation}and a sequence of functions
$\{h_k(s,\theta)\}_{k=0}^\infty$ defined inductively by
\begin{equation}\label{(2.7)}h_{k} =
\alpha_k(\theta)+s\beta_k(\theta)+(s^{2}+1)g_k(s,\theta),~~h_{k+1}
=\frac {s^2+1}\theta\left\{ \frac {\alpha-1}{s+i}+\frac
{\beta-1}{s-i}-\frac d {ds}\right\}g_k(s,\theta)
\end{equation}for $k=0,1,\cdots$. The coefficients $\alpha_k(\theta)$
and $\beta_k(\theta)$ are determined by requiring all
$h_{k}(s,\theta)$ and $g_{k}(s,\theta)$  to be analytic in
$\mathcal{D}=\{s:\Re s\geq-\frac{\eta}{\theta},|s\pm
i|<\frac{2\pi}{\theta}\}$. The uniform asymptotic expansion is given in the following theorem:
\vskip .5cm

\noindent {\thm\label{theorem 1}{For $\theta\in[0,\pi-\delta]$ with
arbitrary  $\delta>0$, we have
\begin{equation}\label{(2.8)}C_{n}^{(\alpha,\beta)}(e^{i\theta})=\theta^{1-\alpha-\beta}T_{1}(n\theta)
\sum_{k=0}^{m-1}\frac{\alpha_{k}(\theta)}{n^{k}}+
\theta^{1-\alpha-\beta}T_{2}(n\theta)\sum_{k=0}^{m-1}\frac{\beta_{k}(\theta)}{n^{k}}
+\varepsilon_{\theta,m},
\end{equation}
where $|\alpha_{k}(\theta)|\leq M_{k}\,,
|\frac{\beta_{k}(\theta)}{\theta}|\leq M_{k}$, for $ k=0,1,2,\cdots
$, and
\begin{equation}\label{(2.9)}|\varepsilon_{\theta,m}|\leq
M_{m} {\theta^{1-\alpha-\beta}} {n^{-m}}\left
\{|T_{1}(n\theta)|+|T_{2}(n\theta)|\right \},
\end{equation}
for $ m=1,2,\cdots $. The positive constants $M_{k}$, $k=1,2,\cdots$,
are independent of  $\theta$ for $\theta\in[0,\pi-\delta]$. The
coefficients $\alpha_{k}(\theta)$ and $\beta_{k}(\theta)$ are given
iteratively in  (\ref{(2.7)}), with $$\alpha_0=\frac
{e^{i\theta\alpha}} 2\left (\frac{\sin\theta}\theta\right
)^{-\alpha}+\frac {e^{-i\theta\beta}} 2\left
(\frac{\sin\theta}\theta\right )^{-\beta},~~ \beta_0=\frac
{e^{i\theta\alpha}} {2i}\left (\frac{\sin\theta}\theta\right
)^{-\alpha}-\frac {e^{-i\theta\beta}} {2i}\left
(\frac{\sin\theta}\theta\right )^{-\beta}.$$}}

\noindent
 \section*{3. Special Functions $T_1$ and $T_2$} \indent\setcounter{section} {3}
\setcounter{equation} {0} \label{sec:3}

 By expanding the slowly varying factor in the integrand of
(\ref{(2.5)}) in uniformly convergent power series of $1/s$  and
integrating   term by
 term, we obtain
 \begin{equation}\label{(3.1)}
                                       T_{1}(x)=
                        x^{\alpha+\beta-1}\sum_{k=0}^{\infty}\sum_{l=0}^{k}\frac{i^k(-1)^{l}(\alpha)_{l}
                        (\beta)_{k-l}}{l!(k-l)!\Gamma(\alpha+\beta+k)}
                        x^k,
                     \end{equation}
where $x^{\alpha+\beta-1}$ is positive for real positive $x$. It is
easily    seen that $x^{-\alpha-\beta+1}T_{1}(x)$ is an entire
function.  From (\ref{(2.5)}), it is also readily verified  that
$T_{2}(x)=T_{1}^{'}(x)$  in the cut plane $x\in \mathbb{C}\backslash
(-\infty, 0]$.

The behavior of    $T_{1}(x)$  at the origin is exhibited  by the leading
terms in (\ref{(3.1)}).  The asymptotic behavior of $T_l(x)$ as
$x\rightarrow +\infty$ has already been given in
\cite[(4.45)]{wongzhao2005}. Moreover, with  $a=2-\alpha-\beta$
and $b=\beta-\alpha$, $T_1(x)$ satisfies the differential
equation
\begin{equation}\label{(3.2)}xT_1''+aT_1'+(x-bi)T_1=0;\end{equation} cf. \cite[(4.52)]{wongzhao2005}.
This  equation is of Laplace type in the sense that the coefficients
are linear in $x$.  We proceed to show that   $T_1(x)$ is connected
to the confluent hypergeometric function via the following change of
variables:
\begin{equation}\label{(3.3)}T_1(x)=  x^{1-a} e^{i x} y(x),~~~z=-2i x.\end{equation} Substituting  (\ref{(3.3)})
in  (\ref{(3.2)}) leads to
  the Kummer equation
\begin{equation}\label{(3.4)} z\frac {d^2 y}{dz^2} +\left [(2-a)-z\right ] \frac {dy}{dz}
-\frac {2-a-b} 2 y=0.\end{equation}   Taking into account the first
two terms of the infinite series in (\ref{(3.1)}), we have
\begin{equation}\label{(3.5)}T_1(x)=
\frac{1}{\Gamma(\alpha+\beta)}x^{\alpha+\beta-1}e^{ix}
M(\alpha,\alpha+\beta,-2ix),\end{equation} where $M$ is the Kummer
function; cf. \cite[p.322]{olver2010}. (One may also derive this
result directly from  (\ref{(2.5)}) and a Laplace inversion integral
of this Kummer function \cite[p.327, (13.4.13)]{olver2010}.)
Accordingly, from $T_{2}(x)=T_{1}^{'}(x)$ one has
\begin{equation}\label{(3.6)}T_2(x)=
\frac{(\alpha+\beta-1+ix)x^{\alpha+\beta-2}e^{ix}
}{\Gamma(\alpha+\beta)} M(\alpha,\alpha+\beta,-2ix) - \frac{2i
x^{\alpha+\beta-1}e^{ix}}{\Gamma(\alpha+\beta)}
M'(\alpha,\alpha+\beta,-2ix) ,\end{equation} where
$M'(\gamma,\delta,z)=\frac d {dz}M(\gamma,\delta,z)$. Thus,
substituting (\ref{(3.5)}) and (\ref{(3.6)}) in (\ref{(2.8)}), we
obtain an asymptotic expansion of the Heisenberg polynomials in
terms of  the Kummer function as follows: \vskip .5cm

\noindent {\thm\label{theorem 2}{ Assume that $\alpha$ and $\beta$
are real and fixed, $z=\rho e^{i\theta}$ with $\rho >0$ and $\theta$
real.  Then we have the compound  asymptotic expansion (\S4.2 of
Chapter 4 in \cite{olver1974})
\begin{equation}\label{(3.7)}C_{n}^{(\alpha,\beta)}(z)\sim n^{\alpha+\beta- 1}z^n \left\{M(\alpha, \alpha+\beta, -2in\theta)
\sum_{k=0}^{\infty}\frac{c_{k}(\theta)}{n^{k}}+
  M'(\alpha, \alpha+\beta, -2in\theta)\sum_{k=0}^\infty\frac{d_{k}(\theta)}{n^{k}}\right\}
  \end{equation}
as $n\rightarrow \infty$, uniformly with respect to $\rho\in (0, \infty)$ and 
$\theta\in[0,\pi-\delta]$, where $\delta$ is an arbitrary constant
such that $0<\delta \leq \pi$. Here the coefficients are given by
\begin{equation}\label{(3.8)}c_k(\theta)=\frac {\alpha_k(\theta) +i \beta_k(\theta)} {\Gamma(\alpha+\beta)}+
\frac { \beta_k(\theta)/\theta }
{\Gamma(\alpha+\beta-1)},~~~d_k(\theta)=\frac {-2i\beta_k(\theta)}
{\Gamma(\alpha+\beta)}
\end{equation}
for $k=0,1,2,\cdots $,   $\alpha_k(\theta)$ and $\beta_k(\theta)$
being defined as in Section 2. In particular,
  $$c_0(\theta)=\frac
{e^{i\theta\alpha}} {\Gamma(\alpha+\beta)}  \left (\frac{\sin\theta}\theta\right
)^{-\alpha},~~~d_0(\theta)=
\frac {e^{-i\theta\beta}} {\Gamma(\alpha+\beta)}\left
(\frac{\sin\theta}\theta\right )^{-\beta}
-
\frac
{e^{i\theta\alpha}} {\Gamma(\alpha+\beta)}\left (\frac{\sin\theta}\theta\right
)^{-\alpha}.$$}}
\vskip .5cm

When the parameters $\alpha$ and $\beta$ in the above theorem are
nonpositive integers, the coefficients $c_k$ and $d_k$ in
(\ref{(3.7)}) all vanish; see (\ref{(3.8)}). However, the asymptotic
relation remains valid since the polynomials also vanish for large
values of $n$; i.e., it is a trivial result.

 For completeness, we write down, by using  (\ref{(2.5)}) and
(\ref{(3.2)}),   the asymptotic expansion  in the special case
$\alpha=\beta$; see also \cite{wongzhao2005}. In this case, we have
\begin{equation}\label{(3.9)}
T_1(x)=\frac{\sqrt\pi}{\Gamma(\alpha)} \left (\frac x 2 \right )^{\alpha-  1 /2} J_{\alpha- 1/ 2} (x),~~~
 T_2(x)=\frac{\sqrt\pi}{\Gamma(\alpha)} \left (\frac x 2 \right )^{\alpha-  1 /2} J_{\alpha-  3 /2} (x),
\end{equation}
where $J_\nu(x)$ is the Bessel function of the first kind; cf.
\cite[p.217]{olver2010}. The expansion in Theorems \ref{theorem 1}
and  \ref{theorem 2} now takes the form: \vskip .5cm

\noindent {\cor\label{corollary 1}{ Assume that $\alpha$ is  real
and fixed, $z=\rho e^{i\theta}$ with $\rho >0$ and $\theta$ real.
Then we have the compound asymptotic expansion
\begin{equation}\label{(3.10)}C_{n}^{(\alpha,\alpha)}(z)\sim \rho^n
\frac{\sqrt\pi}{\Gamma(\alpha)}
\left (\frac n {2\theta} \right )^{\alpha-1/2}
 \left\{J_{\alpha- 1/ 2} (n\theta)
\sum_{k=0}^{\infty}\frac{\alpha_{k}(\theta)}{n^{k}}+
  J_{\alpha- 3/ 2} (n\theta)\sum_{k=0}^\infty\frac{\beta_{k}(\theta)}{n^{k}}\right\} 
\end{equation}
as $n\rightarrow \infty$, uniformly with respect to $\rho\in (0, \infty)$ and 
$\theta\in[0,\pi-\delta]$, where $\delta$ is an arbitrary constant
such that $0<\delta\leq \pi$.   Here $\alpha_k(\theta)$ and
$\beta_k(\theta)$ are given in (\ref{(2.7)}), and the leading
coefficients are
\begin{equation}\label{(3.11)}\alpha_0(\theta)=\cos(\alpha \theta)\left (\frac {\sin\theta}\theta\right )^{-\alpha}~~~
\mbox{and}~~~\beta_0(\theta)=\sin(\alpha \theta)\left (\frac {\sin\theta}\theta\right )^{-\alpha}.\end{equation}}}
\vskip .5cm

\noindent
 \section*{4. Asymptotic Zeros} \indent\setcounter{section} {4}
\setcounter{equation} {0} \label{sec:4}

We complete this note by giving some remarks on the  zeros
of the Heisenberg polynomials. By (\ref{(2.8)}),
these zeros
 are approximated by the zeros
of $T_1(n\theta)$. Thus,  we shall later  turn to consider the zeros of $T_1$.  But, as a consistency check on our result, we first consider the special
case $\alpha=\beta >- 1/2$.   In this case, the Heisenberg polynomials on the unit circle are related to the Gegenbauer
polynomials via the identity
\begin{equation}\label{(4.1)}C_{n}^{(\alpha,\alpha)}(e^{i\theta})=C_{n}^{(\alpha)}(\cos\theta), ~~\theta \in [0, \pi],
\end{equation}
which  can be seen by comparing the generating functions in
(\ref{(1.6)}) and in \cite[p.448]{olver2010}. In this specific case,
$C_{n}^{(\alpha,\alpha)}(e^{i\theta})$ is known to have $n$ distinct
  zeros in the interval $0< \theta <\pi$.

The asymptotic behavior of the extreme zeros can be obtained from  (\ref{(3.10)}) and (\ref{(3.11)}). Indeed, when $n\theta$ is bounded,  the terms in the curly brackets in   (\ref{(3.10)}) can be written as
$H:=\left (\alpha_0+ {\alpha_1}/ n\right )J_{\alpha- 1/ 2} (n\theta)
 +\beta_0
  J_{\alpha- 3/ 2} (n\theta)+O(n^{-2})$; here,  we have used the fact that $\beta_k(\theta)=O(\theta)=O(1/n)$ for $k=0,1,2,\cdots$,
   as mentioned in Theorem \ref{theorem 1}. The quantity $H$  can be further  expressed  as
  \begin{equation}\label{(4.2)}H=\Omega \left [ J_{\alpha-1/2}(n\theta)+\Delta
J'_{\alpha-1/2}(n\theta)+O(n^{-2})\right ]=\Omega \left [ J_{\alpha-1/2}(n\theta+\Delta)
+O(n^{-2})\right ],
\end{equation}
where $\Omega=  \alpha_0+\{\alpha_1+ (\alpha-\frac 1 2) {\beta_0}/\theta \}n^{-1}$  and $\Delta=\beta_0/\Omega$  both  depend  on
$\theta$ and $n$. Note  that $\Omega=\alpha_0+O(1/n)=1+O(1/n)$ and $\Delta=\left\{ n\theta \Omega^{-1}  \beta_0/\theta \right\} / n=O(1/n)$ as  $n\rightarrow\infty$. In  (\ref{(4.2)}), we have used the Taylor expansion of
$J_{\alpha-1/2}(x)$ at $x=n\theta$   and  the  formula
$$J_{\alpha-3/2}(x)=J'_{\alpha-1/2}(x)+\frac {2\alpha-1} {2x} J_{\alpha-1/2}(x);$$
cf. \cite[p.222]{olver2010}. Hence, $C_{n}^{(\alpha,\alpha)}(z)=0$
amounts to $H=0$, which implies
$$J_{\alpha-1/2}(n\theta+\Delta)
=O(n^{-2}).$$
From  (\ref{(3.11)}),
it is readily seen  that $\Delta=\alpha\theta+O(1/n^2)$. Denoting by $j_{\alpha- 1/ 2, k}$ the
$k$th positive zero of $J_{\alpha- 1/ 2}(x)$, we have the asymptotic approximation
\begin{equation}\label{(4.3)}\theta_{n,k}=\frac {j_{\alpha-  1 /2, k}}{n+\alpha} +O(n^{-3})
\end{equation}
for the $k$th zero of $C_{n}^{(\alpha,\alpha)}(e^{i\theta})$  arranged in increasing order, where $k=1,2,\cdots$ are fixed and  $\alpha>-1/2$. In deriving  (\ref{(4.3)}), we have made use of the fact that $J'_{\alpha-1/2}$ does not vanish in a neighborhood of $j_{\alpha-  1 /2, k}$.

If  we relax  the restriction that $k$ be fixed,  while still requiring  $\theta=O(1/n)$ and  $n\gg k \gg 1$, then formula (\ref{(4.3)}) would continue to  hold, with only a weaker error estimate. For instance, making  use of the fact that
$$j_{\alpha- 1 /2, k}\sim  \left ( {k+(\alpha-1)/2}\right )\pi $$
for $\alpha>1/2$ and $k\gg 1$ (cf. \cite[p.236,
(10.21.19)]{olver2010}), and substituting  this into (\ref{(4.3)}),
gives
\begin{equation}\label{(4.4)} \theta_{n,k}\sim \frac { (\alpha-1)/ 2  +k} {n+\alpha} \pi
\end{equation}for   $\alpha>1/2$ and $n\gg k \gg 1$. The formula in (\ref{(4.4)}) agrees with that in
\cite[p.139, (6.6.3)]{szego1975} for $1/2<\alpha<1$.

On the other hand, taking   the two  leading terms in
(\ref{(3.10)})-(\ref{(3.11)}) and making use of the asymptotic
approximation for the Bessel function \cite[p.223,
(10.7.8)]{olver2010}
\begin{equation}\label{(4.5)} J_\nu(x)=\sqrt{2/(\pi x)}\cos(x- \nu\pi/2- \pi/4) +O(x^{-1})~~\mbox{as}~~x\rightarrow +\infty,
\end{equation}
we have  from (\ref{(3.10)})
\begin{equation}\label{(4.6)}C_{n}^{(\alpha,\alpha)}(z)= \rho^n
\frac{\theta^{-\alpha}}{\Gamma(\alpha)}
\left (\frac n 2 \right )^{\alpha-1}\left (\frac {\sin\theta}\theta\right )^{-\alpha}
 \left\{\cos\left ((n+\alpha)\theta-\frac {\alpha\pi} 2\right )+O\left (\frac 1 {n\theta}\right ) \right\}
\end{equation}
for $n\theta\gg 1$ and $\theta \in (0, \pi-\delta]$. It is readily seen that this asymptotic approximation
agrees well   with (\ref{(4.4)}).

It is worth mentioning that if $\theta=\theta_{n,\nu}$ is a zero of $C_{n}^{(\alpha,\alpha)}(e^{i\theta})$, then all points on the ray $\arg z=\theta_{n,\nu}$ are zeros of $C_{n}^{(\alpha,\alpha)}(z)$. This fact follows readily from the homogeneity of the
Heisenberg polynomials, namely,   $C_n^{(\alpha, \beta)}(\rho e^{i\theta} )=\rho^n C_n^{(\alpha, \beta)}( e^{i\theta} )$ for $\rho>0$ and
$\theta\in (-\pi, \pi]$.

Now,  we turn to the asymptotic zeros of the special function $T_1(x)$,   whose asymptotic properties
were investigated   in \cite{wongzhao2005}. It is of
interest to consider the zeros of $T_1(x)$ for $x>0$ when
$\alpha\neq\beta$. Without loss of generality,  we   assume that $\beta
>\alpha$.

{\pro\label{proposition 1}{The function $T_1(x)$ possesses  no
positive zeros   when $\alpha+\beta>0$. }}\\

\noindent {\bf{Proof:}}  We use an argument similar to that given in Ince
\cite[p.512]{ince1956}. Write
$T_1(x)=x^{\frac 1 2 (\alpha+\beta)-1}I(x)$. Equation (\ref{(3.2)}) gives
\begin{equation}\label{(4.7)}I^{''}(x)+[g_{1}(x)+ig_{2}(x)]I(x)=0,\end{equation}
where the functions
$$ g_{1}(x)=1-\frac{a(a-2)}{4x^2}~~~ \mbox{and}~~~ g_{2}(x)=-\frac b x$$
are real, and $g_2(x)< 0$ for $x>0$. From (\ref{(4.7)}) we have
$$\left (\overline{I} I'\right )'= \left | I'\right |^2 + \overline
{I} I''=\left \{\left | I'\right |^2 - g_1|I|^2 \right\} -i\left\{
g_2(x) |I|^2\right\}.$$ Integrating over an arbitrary interval
$[x_1, x_2]$,  and taking the imaginary part, one obtains
$$ \Im \left [ \overline {I(x)} I'(x)\right ]
_{x_1}^{x_2} +\int^{x_2}_{x_1} g_2(x) |I(x)|^2 dx=0.$$
 From (\ref{(3.1)}), it is
straightforward to verify that
when
$\alpha+\beta>0$,
$\lim_{x\rightarrow 0^+} \Im\left [ \overline {I(x)} I'(x)\right
]=0$. Hence, for $x>0$,
$$\Im \left [ \overline {I(x)} I'(x)\right ]  =  - \int^{x}_{0} g_2(s) |I(s)|^2
ds>0,$$ from which it follows that  $I(x)$, and consequently
$T_1(x)$, has no zero in the interval $(0, \infty)$. \vskip .5cm

In view of Theorem  \ref{theorem 2},   we conclude that for large values of $n$ there are no  zeros of
$C_{n}^{(\alpha,\beta)}(z)$ in the complex $z$-plane when $\alpha\not=\beta$ and $\alpha+\beta>0$. In our discussion of asymptotic zeros of
$C_{n}^{(\alpha,\beta)}(z)$,
we have not addressed the cases when $\alpha=\beta\leq -1/2$  and  when $\alpha\not=\beta$ with $\alpha+\beta \leq 0$.

\end{document}